\documentstyle[amstex]{article}

\begin{document}

\title{{\bf   Realcompactness and Banach-Stone theorems}}

\author{\normalsize {\sc Jes\'us Araujo}\thanks{Research partially 
supported by the
Spanish Direcci\'on General de Investigaci\'on Cient\'{\i}fica y
T\'ecnica (DGICYT,  PB98-1102).}}

\date{}                                         

\maketitle 








\newtheorem{theorem}{Theorem}[section]
\newtheorem{corollary}[theorem]{Corollary}
\newtheorem{lemma}[theorem]{Lemma}
\newtheorem{proposition}[theorem]{Proposition}
\newtheorem{ax}{Axiom}

\newtheorem{definition}{Definition}[section]

\numberwithin{equation}{section}

\newtheorem{examples}{Examples}

\newtheorem{example}{Example}
















\newcommand{\cly}{{\rm cl}_{\beta Y} \hspace{.02in}}
\newcommand{\cgy}{{\rm cl}_{\gamma Y} \hspace{.02in}}
\newcommand{\chy}{{\rm cl}_{Y} \hspace{.02in}}
\newcommand{\clx}{{\rm cl}_{\beta X} \hspace{.02in}}
\newcommand{\cgx}{{\rm cl}_{\gamma X} \hspace{.02in}}
\newcommand{\clp}{{\rm cl}_{{\Bbb R}^p} \hspace{.02in}}
\newcommand{\chx}{{\rm cl}_{X} \hspace{.02in}}
\newcommand{\cux}{{\rm cl}_{\upsilon X} \hspace{.02in}}
\newcommand{\bx}{A(X,E)}
\newcommand{\xe}{C_b(X,E)}
\newcommand{\yf}{C_b(Y,F)}
\newcommand{\ra}{\rightarrow}
\newcommand{\smn}{\sum_{n=1}}
\newcommand{\ay}{A(Y,F)}
\newcommand{\cxo}{C_0(X,E)}
\newcommand{\cyo}{C_0(Y,F)}
\newcommand{\cx}{C_0(X)}
\newcommand{\sxy}{T:A(X) \ra A(Y)}
\newcommand{\txy}{T:A(X,E) \ra A(Y,F)}
\newcommand{\bxy}{T:C_b (X,E) \ra C_b (Y,F)}
\newcommand{\uxy}{T:C(X,E) \ra C(Y,F)}
\newcommand{\vxy}{T:C_0 (X,E) \ra C_0 (Y,F)}
\newcommand{\pa}{\left(}
\newcommand{\rb}{\right)}
\newcommand{\va}{\left|}
\newcommand{\vb}{\right|}
\newcommand{\vc}{\left\|}
\newcommand{\vd}{\right\|}
\newcommand{\cl}{{\rm cl} \hspace{.02in}}
\newcommand{\itr}{{\rm int} \hspace{.02in}}
\newcommand{\itrx}{{\rm int}_{\beta X} \hspace{.02in}}
\newcommand{\bax}{^{\beta X}}
\newcommand{\bay}{^{\beta Y}}
\newcommand{\ing}{{\rm ifn} \hspace{.02in}}
\newcommand{\fin}{{\rm fin} \hspace{.02in}}
\newcommand{\bdry}{{\rm bdry} \hspace{.02in}}
\newcommand{\hs}{\hspace{.02in}}

\newcommand{\noe}{A^n (\Omega, E)}
\newcommand{\coe}{C^n (\Omega, E)}
\newcommand{\no}{C^n (\Omega, {\Bbb  R})}
\newcommand{\nor}{C^n ({\Bbb  R}, {\Bbb  R})}
\newcommand{\noor}{C^n_c (\Omega, {\Bbb  R})}
\newcommand{\mor}{C^m ({\Bbb  R}, {\Bbb  R})}
\newcommand{\moor}{C^m_c (\Omega', {\Bbb  R})}
\newcommand{\mo}{C^m (\Omega', {\Bbb  R})}
\newcommand{\moe}{A^m (\Omega', F)}
\newcommand{\bnoe}{C^n_* (\Omega, E)}
\newcommand{\bno}{A^n (\Omega, {\Bbb  R})}
\newcommand{\bmo}{C^m_* (\Omega', {\Bbb  R})}
\newcommand{\bmoe}{C^m_* (\Omega', F)}
\newcommand{\ixy}{T: \bnoe \ra \bmoe}
\newcommand{\partl}{\frac{\partial^{\va \lambda \vb} l}{\partial x^{\lambda}}}
\newcommand{\parto}{\frac{\partial^{\va \lambda \vb} l'}{\partial x^{\lambda}}}
\newcommand{\partm}{\frac{\partial^{\va \lambda \vb} l}{\partial x^{\lambda}}}
\newcommand{\partg}{\frac{\partial^{\va \lambda \vb} Tg}{\partial x^{\lambda}}}
\newcommand{\partf}{\frac{\partial^{\va \lambda \vb} Tf}{\partial x^{\lambda}}}





\date{}                                          
\thispagestyle{empty}

\maketitle

\begin{abstract}
For realcompact spaces $X$ and $Y$
we give a complete description of the linear biseparating maps between spaces of vector-valued continuous functions on $X$ and $Y$, where special attention is paid to spaces of vector-valued {\em bounded} continuous functions. 
These results are applied to describe  the linear isometries between spaces of vector-valued bounded continuous and uniformly continuous functions.
\end{abstract}

\section{Introduction}
Let ${\Bbb K} = {\Bbb R}$ or ${\Bbb C}$.
Given a completely regular space $X$, and a ${\Bbb K}$-normed space $E$,
$C(X,E)$ and $\xe$ denote the spaces of continuous
functions and {\em
bounded} continuous functions on $X$ taking 
values on $E$, respectively. $C(X)$ and $C_b (X) $ will be the
spaces $C(X, {\Bbb K})$ and $C_b (X, {\Bbb K})$, respectively.

(Bi)separating linear maps between spaces of scalar-valued continuous functions have drawn attention of researchers recently. In general, they can be described as weighted composition maps (see for instance \cite{AK}, \cite{A1}, \cite{ABN1}, \cite{ABN2} and \cite{J}). As a result, automatic continuity for this kind of maps is obtained as a corollary.

As for spaces of vector-valued continuous functions, a similar approach is taken for biseparating linear maps in \cite{HBN} and \cite{A2}, where a description as weighted composition maps is obtained when topological spaces $X$ and $Y$ are compact or locally compact.

In this paper, we drop every assumption of (local) compactness on $X$ and $Y$, assuming realcompactness instead. Notice at this point that the class of realcompact spaces is fairly large since it includes, apart from that of compact spaces, the class of subsets of Euclidean spaces and even the class of all metric spaces of nonmeasurable cardinal (see for instance \cite[p. 232]{GJ}). When dealing with realcompact spaces we obtain a representation of biseparating linear maps similar to that given in the compact setting. Special mention deserves the fact that this description of biseparating linear maps as weighted compositions apply even when they are just defined between spaces of {\em bounded} continuous functions $\xe$ and $\yf$, whenever ·$E$ and $F$ are infinite-dimensional (see Theorem~\ref{linear}). 

        The other aim of this paper is the study of a topological link (namely, a linear isometry) between spaces of bounded continuous functions yielding also a topological relation between $X$ and $Y$, when they are realcompact.  A classical result in this theory turns out to be
the following, which allows us to describe the surjective linear
isometries
$T: C(X) \ra C(Y)$, where $X$ and $Y$ are assumed to be compact:
if $T:C(X) \ra C(Y)$ is a surjective linear isometry, then there
exist a homeomorphism $h$ from $Y$ onto $X$ and $a \in C(Y)$
with $\va a(y) \vb=1$ for every $y \in Y$ such that\[ (Tf) (y) =
a(y) f(h(y))\] for each $f \in C(X)$ and each $y \in Y$.

On the other hand,  a linear isometry $T:C_b(X) \ra C_b(Y)$ always
lead to a homeomorphism between $\beta X$ and $\beta Y$ (the Stone-\v{C}ech compactifications of $X$ and $Y$), as elements in $C_b(X)$ and $C_b(Y)$ can be extended to elements in $C(\beta X)$ and $C(\beta Y)$. 
For a
systematic account on isometries between spaces of continuous functions, see for instance \cite{FJ} and \cite{JP2}.

        Also for compact $X$ and $Y$, a decisive step forthward was
taken by Behrends  (see for instance \cite{B}). It says in particular that
whenever $E$ satisfies a special condition, then every
linear isometry from $C(X,E)$ onto $C(Y,E)$ is a {\em strong
Banach-Stone map} (see definition in Section 3), that is, in that case
we can obtain a description of the map. This result could be extended
to maps defined between spaces of continuous functions vanishing
at infinity on
locally compact spaces. In general, results in this trend always
include some kind of compacity of the topological spaces among the
hypotheses. 

        As for isometries between spaces
of vector-valued bounded continuous functions defined on
nonlocally compact spaces, we would like to point out that the only results which
seem to have made its way in the literature so far are contained in \cite{Bachir}, where the author gives a representation of such isometries in the case when $X$ and $Y$ are complete metric spaces and $E$ is a Hilbert space. 

Our approach here is different:  we take advantage of our study of biseparating maps to describe such isometries. But not only this, because a similar approach can be taken to study linear isometries between spaces of bounded uniformly continuous functions, providing in this case a special description for them (see Theorem~\ref{gol}).

        Summing up, in this paper we drop the assumption of compacity or
local compacity on $X$
and $Y$, and arrive at the same conclusion: for some kind of
Banach spaces $E$ and $F$, every surjective
linear isometry $\bxy$ is a strong Banach-Stone map, and when the isometry is defined between spaces of bounded uniformly continuous functions, then both the homeomorphism between $X$ and $Y$ and its inverse are uniformly continuous.

	On the other hand, we also mention that the requirements of realcompactness on our spaces is necessary for the descriptions we provide. If $X$ or $Y$ are not realcompact, in general the biparating linear maps or the linear isometries do not admit such representations. An example can be given even in a very easy context:

\bigskip

\noindent
{\bf Example.}
Assume that $X$ is not realcompact (for instance $X= W( \omega_1):= \{\sigma : \sigma < \omega_1\}$, where $\omega_1$ denotes the first uncountable ordinal; see \cite[5.12]{GJ}), and $E= l^2$. Since $l^2$ is separable, it is realcompact (see \cite[8.2]{GJ}). Consequently each bounded continuous map $f: X \ra l^2$ can be extended to a continuous map $f^{\upsilon X} : \upsilon X \ra l^2$ defined in the realcompactification $\upsilon X$ of $X$, which is also bounded. In this way we can define a linear isometry (which is also a biseparating map) from $C_b (X, E)$ onto $C_b (\upsilon X, E)$. Clearly this map does not admit a description as the one given in Theorems~\ref{linear} and \ref{gol}.

\bigskip

In the rest of this section and in Section 2 $E$ and $F$ will be ${\Bbb K}$-normed spaces, whereas in Section 3 they will be ${\Bbb K}$-Banach spaces. All over the paper $X$ and $Y$ will be completely regular topological spaces. Finally, if $X$ and $Y$ are also complete metric spaces, we introduce  $C_b^u(X, E)$ and $C_b^u(X, E)$ as the spaces  of uniformly continuous bounded functions defined on $X$ and $Y$, respectively, and taking values in $E$ and $F$, respectively. Also in this case $C_b^u(X) = C_b^u(X, {\Bbb K})$, and $C_b^u(Y) = C_b^u(Y, {\Bbb K})$.

\begin{definition}
{\em Given $f \in C(X,E)$, we define the cozero set of $f$ as \[ c(f) := \{ x \in X : f(x) \neq 0\}.\]}
\end{definition}


\begin{definition}
{\em Let $A(X,E) \subset C(X,E)$, $A(Y,F) \subset C(Y,F)$. A map $\txy$ is said to be {\em separating} if it is additive and
$c(Tf) \cap c(Tg) = \emptyset$ whenever $f, g \in \bx$ satisfy
$c(f) \cap c(g) = \emptyset$. Besides $T$ is said to be
{\em biseparating} if it is bijective and both $T$ and $T^{-1}$ are separating.} 
\end{definition}

Equivalently, we see that an additive map $\txy$ is separating if $\vc (Tf) (y) \vd \vc (Tg) (y) \vd =0$ for all $y \in Y$ whenever $f, g \in \bx$ satisfy $\vc f(x) \vd \vc g(x) \vd=0$ for all $x \in X$. 

       

As for the spaces of linear functions, we will denote by $L' (E,F)$ and by $B' (E,F)$ the sets of (not necessarily continuous) linear maps and  bijective linear maps  from $E$ into $F$, respectively.  On the other hand, $B(E,F)$ and  $I(E,F)$ stand for the set of continuous bijective linear maps from $E$ into $F$ and the set of all linear isometries from $E$ onto $F$, respectively. We consider that both $B(E,F)$ and $I(E,F)$ are endowed with the strong operator topology, that is, the coarsest topology such that the mappings $J \hookrightarrow J {\bf e}$ are continuous for every ${\bf e} \in E$.

	All over the paper the word "homeomorphism" will be synonymous with "surjective homeomorphism".

\section{Representation of linear  biseparating maps}

        In this section we assume that $E$ and $F$ are ${\Bbb K}$-normed spaces, and that we are in one of the following situations:

\begin{itemize}
\item {\bf Situation  1.} $\bx = C(X,E)$ and $\ay = C(Y,F)$. Also $X$ and $Y$ are asumed to be realcompact.
\item {\bf Situation  2.} $E$ and $F$ are infinite-dimensional, and $\bx= C_b (X,E)$, $\ay = C_b (Y,F)$. As above, also in this case $X$ and $Y$ are realcompact.
\item {\bf Situation  3.} $\bx = C_b^u(X, E)$ and $\ay = C_b^u(Y, F)$. In this case $X$ and $Y$ are complete metric spaces.

\end{itemize}

We also assume that if we are in Situations 1 or 2, then $A := C_b (X)$. Otherwise  $A:= C_b^u(X)$.

\medskip

In a much more general context, (not necessarily linear) biseparating maps are studied in \cite{A3}. Among other things, it is proven there that the existence of a biseparating map from $A(X,E)$ onto $A(Y,F)$ leads to the existence of a homeomorphism $h: Y \ra X$. Among the properties of this map $h$ (called {\em support map}) we have the following:

\begin{lemma}\label{13}(\cite[Lemma 4.4]{A3})
Let $\txy$ be a  biseparating map. Suppose that $h(y)
= x$ for some $y \in Y$, and that $f \in \bx$
satisfies $f  \equiv 0$ on a neighborhood of
$x$. Then $Tf \equiv 0$ on a neighborhood of $y$.
\end{lemma}

\begin{corollary}\label{unif}(\cite[Corollary 3.3]{A3})
Let $\bx = C_b^u(X, E)$, $\ay = C_b^u(Y, F)$. If $\txy$ is a biseparating map, then $h: Y \ra X$ is a uniform homeomorphism (that is, both $h$ and $h^{-1}$ are uniform maps).
\end{corollary}

Thanks to these basic results, we can prove the following proposition.

\begin{proposition}\label{14}
Suppose that $\txy$ is  biseparating and that \[f
 (h(y))
=0\] for some $f \in \ay$ and $y \in Y$. Then $(Tf)(y) =0$.
\end{proposition}

{\em Proof.} 
By Lemma~\ref{13}, the result holds if $ f \equiv 0$ on
a neighborhood of $h(y)$. Thus we can assume that $ f$ is
not constant in any neighborhood of $h(y)$. Then consider the
following sequence of neighborhoods of $h(y)$. For $n \in
{\Bbb N}$, let \[U_n := \left\{x \in X :  \vc f
(x) \vd <
\frac{1}{n^3}\right\}.\] 
Without loss of generality we may assume that $U_n \neq U_m$
whenever $n, m \in {\Bbb N}$, $n \neq m$. It is also clear that
\[h(y) \in V_1:= \cl \bigcup_{n\in{\Bbb N}}  (U_{4n} -
\cl  U_{4n+3})\] or \[h(y) 
\in V_2:= \cl \bigcup_{n \in {\Bbb N}}  (U_{4n-3} -  \cl U_{4n}).\]

        We assume without loss of generality that $h(y)$ belongs to
$V_1$. Notice that in this case \[ h(y) \in \cl
\bigcup_{n \ge k}  (U_{4n} - 
\cl  U_{4n+3})\] for every $k \in {\Bbb N}$.

Consider a sequence $(f_n)$ in $A$ such that, for every
$n \in {\Bbb N}$, $0 \le f_n \le 1$, $c(f_n) \subset U_{4n-1}$,
and $f_n \equiv 1$
in $U_{4n}$. It is clear that $g:= \sum_{n=1}^{\infty} f_n f$
belongs to $\bx$ because $\sup_{x \in X} \vc f_n f (x) \vd \le
1/n^2$ and $g(x) \in E$ for every $x \in X$. On the other hand, it is easy to see that $g \equiv nf $
on $U_{4n} - \cl U_{4n+3}$ for each $n \in {\Bbb N}$.

        Next suppose that $(Tf) (y) = {\bf e}_0 \in F$, ${\bf e}_0 \neq
0$,  and $(Tg) (y) = {\bf e}_1 \in F$. 
Consider $n_0 \in
{\Bbb N}$ with \[n_0 \vc {\bf e}_0 \vd /2 > \vc {\bf e}_1 \vd +1,\]
and an open
neighborhood $U(y)$ of $y$ in $Y$ such that $h(U(y)) \subset
U_{4n_0}$ and \[\vc (Tf) (y') \vd > \vc {\bf e}_0 \vd /2\]
for every $y' \in U(y)$.


 Taking into account that $h$ is a homeomorphism and that $h(y)$ belongs to 
$\cl
\bigcup_{n \ge k}  (U_{4n} - 
\cl  U_{4n+3})$ for every $k \in {\Bbb N}$,
we can see that there exists $k \in
{\Bbb N}$, $k \ge n_0$, such that \[h(U(y)) \cap (U_{4k} -
\cl
U_{4k-3}) \neq \emptyset.\]

 Then if for  $y_1 \in U(y)$, $h(y_1)$ belongs
to $U_{4k} - \cl 
U_{4k-3}$, we have that $ g  - kf$ is constantly
equal to zero in
a neighborhood of $h(y_1)$. We deduce by Lemma~\ref{13} that
$Tg \equiv k Tf$ on a neighborhood of $y_1$,
which implies that \[\vc (T g)
 (y_1) \vd = k \vc (Tf)  (y_1) \vd > \vc {\bf e}_1
\vd +1.\]

        This behaviour in every neighborhood of $y$ implies that
$Tg$ is not continuous, which is not
possible. We conclude that $(Tf)(y) =0$.
\hfill $\Box$

\begin{theorem}\label{linear}
Suppose that $\txy$ is a linear  biseparating map. Then there
exist a  homeomorphism $h: Y \ra X$ and a map $J: Y \ra B' (E,F)$
such that $(Tf) (y) = (Jy) (f(h(y))$ for every $f \in \bx$ and
$y \in Y$. On the other hand, if we are in Situation 3, then $h$ is also a uniform homeomorphism.
\end{theorem}

{\em Proof.} 
For each $y \in Y$, we define a linear map $Jy :E \ra F$ as
$(Jy) ({\bf e}) = (T \widehat{{\bf e}}) (y)$. It is clear that, if $h:Y \ra X$ is the
support map, then for every $y \in Y$ and
$f \in \bx$, $f(h(y)) = \widehat{f(h(y))} (h (y))$, and by
Proposition~\ref{14}, $(Tf)(y) = (T \widehat{f(h(y))} ) (y) $, that
is, \[ (Tf) (y) = (Jy)
(f(h(y))).\]

Next we prove that each $Jy : E \ra F$ is bijective. Notice
first that
since $T^{-1}$ is also  biseparating, the above representation
can be applied to $T^{-1}$. This 
implies in particular that there exists $K: X \ra L'(F,E)$ such that,
for every $g \in \ay$ and $x \in 
X$, \[(T^{-1} g) (x) = (Kx) (g(h^{-1} (x)).\]

        Fix $y \in Y$ and ${\bf f} \in F-\{0\}$. Let $x =
h(y)$. Now take $g \in 
\ay$ with $g(y) = {\bf f}$. Then its is clear that ${\bf f}
=g (y) =  (T (T^{-1} g)) (y)$, that is,
\begin{eqnarray*}\label{reves}
{\bf f} &=& (Jy) ((T^{-1} g) (x))\\
&=& (Jy) ((Kx) (g(h^{-1} (x))))\\
&=& (Jy) ((Kx) (g(y)))\\
&=& (Jy) ((Kx) ({\bf f})).
\end{eqnarray*}

        This implies that $(Jy) (Kx)$ is the identity map on
$F$. In the same way we can prove that $(Kx)(Jy)$ is the
identity map on $E$. Consequently, $Jy$ is bijective.

Finally, if we are in Situation 3, the fact that $h$ is a uniform homeomorphism follows from Corollary~\ref{unif}.
\hfill $\Box$ 

\bigskip

Next we are going to see that, when we deal with finite-dimensional $E$, some properties regarding continuity can be obtained. The following result follows immediately from Theorem~\ref{linear}.

\begin{corollary}\label{15}
Suppose that $\txy$ is  biseparating and linear, and that $E$ is
finite-dimensional. Then $E=F$.
\end{corollary}

As a consequence of  the last corollary, when $E$ is finite-dimensional, the map $J: Y \ra B'(E,F)$ given in Theorem~\ref{linear} attains values in $B(E,F)$.

\begin{corollary}\label{16}
Suppose that $\txy$ is  biseparating and linear, and that $E$ is
finite-dimensional. Then the map $J: Y \ra B (E,F)$  is continuous.
\end{corollary}

{\em Proof.} 
We have that $E$ is finite-dimensional and, by Corollary~\ref{15}, so
is $F$. Suppose that $(y_i)_{i \in I} $ is a net
in $Y$ which converges to $y_0 \in Y$. To prove that $J$ is
continuous it is enough to prove that, for every ${\bf e} \in E$,
$((Jy_i) ({\bf e}))_{i \in I}$ converges to $(Jy_0) ({\bf e})$. This is clear from the definition of $J$.
\hfill $\Box$

\begin{corollary}\label{200}
Suppose that $\uxy$ is  biseparating and linear, and that $E$ is
finite-dimensional. If $C(X,E)$
and $C(Y,F)$ are endowed with the compact-open topology, then $T$ (and $T^{-1}$)
is continuous.
\end{corollary}

{\em Proof.} 
Let $K$ be a compact subset of $Y$. Recall that, by Theorem~\ref{linear}, $Jy$ is bijective for each $y \in Y$, and that, by
Corollary~\ref{16}, 
 the map $Jy :Y \ra B(E,F)$ is continuous. Suppose that $\epsilon > 0$
and that for every $x \in h(K)$, \[ \vc f(x) \vd <
\frac{\epsilon}{\sup_{y \in K} \vc J y \vd}.\] Now it is
straightforward to see
that $\vc (Tf) (y) \vd < \epsilon$.
\hfill $\Box$ 

\bigskip

The following corollary follows easily from the principle of uniform boundedness.

\begin{corollary}\label{202}
Suppose that $T: C_b (X,E) \ra C_b(Y,F)$ (respectively, $T: C_b^u(X,E) \ra C_b^u(Y,F)$)  is  biseparating and linear, and that $E$ is
finite-dimensional. If $C_b(X,E)$
and $C_b(Y,F)$ (respectively, $C_b^u(X,E) $ and $C_b^u(Y,F) $), are endowed with the sup norm, then $T$
(and $T^{-1}$) is continuous.
\end{corollary}

But in the case of continuous linear  biseparating maps, we can say even more.

\begin{corollary}\label{cont}
Suppose that $\bxy$ (respectively, $T: C_b^u(X,E) \ra C_b^u (Y,F)$)   is a continuous linear  biseparating map. Then the map $J$ given in Theorem~\ref{linear} takes values in $B(E,F)$. Also $J: Y \ra B (E,F)$ is continuous.
\end{corollary}

{\em Proof.} 
We first have to prove that $Jy $
is continuous for each $y \in Y$. But this is clear
because given any ${\bf e} \in E$ and $y \in Y$, \[\vc (Jy)
({\bf e}) \vd = \vc (T \widehat{{\bf e}}) (y)\vd \le \vc T \vd
\vc {\bf e} \vd.\]

	On the other hand,  the map $J: Y \ra B(E,F)$ is continuous, where as we assumed in the Introduction, $B(E,F)$ is endowed with the strong operator topology. 
This can be proved as in Corollary~\ref{16}.
\hfill $\Box$ 

\bigskip

{\bf Remarks.}

        1. When $E$ is
infinite-dimensional, a linear biseparating map need not be continuous,
as it is easy to see in the following example. Consider $E:=
c_0$, the space of sequences converging to zero, and
$X=\{x\}=Y$. Take ${\frak U}$ a Hamel base of $c_0$ such
that every element of ${\frak U}$ has norm one. Consider
${\frak V}:= \{u_n : n \in {\Bbb N}\}$ a countable subset of
${\frak U}$ and define $T u:= u$ if $u \in {\frak U} -
{\frak V}$ and $Tu_n = nu_n$ for $u_n \in {\frak V}$. It is
clear that, by identifying each map in $C( \{x\} , c_0)$ with
its image, $T$ can be extended by linearity to a (clearly
biseparating) bijective map $T:C(\{x\} , c_0) \ra C(\{x\}, c_0)$
which is not continuous.

        2. The existence of a nonlinear biseparating map $T:
C(X,E) \ra C(Y, F)$ does not imply in general that  $E$ and $F$ are isomorphic as vector spaces, even if they are
finite-dimensional. Consider for instance $E := {\Bbb K}$ and $F=
{\Bbb K}^2$. Take a Hamel base ${\frak U} = \{ a_i : i \in I\}$ of
$E$ as a ${\Bbb Q}$-linear space, where ${\Bbb Q}$ is the field of
rational numbers. Clearly ${\frak V} := \{(a_i, 0) : i \in I
\} \cup \{(0, a_j ) : j \in I\}$ is a Hamel base of $F$ as a
${\Bbb Q}$-linear space. Also it
is easy to see that ${\frak U}$ and ${\frak V}$ have the same
cardinal and there exists a bijective map $v : {\frak U} \ra
{\frak V}$. Then we can extend $v$ by ${\Bbb Q}$-linearity to a
bijection defined in the
whole space $E$. Now suppose that $X = \{x\} =Y$. We clearly have
that $T: C(X, E) \ra C(Y, F)$ defined as $T \widehat{{\bf e} } :=
\widehat{v({\bf  e})}$, for every ${\bf e}  \in E$, is a biseparating map
(which obviously is not ${\Bbb K}$-linear).

        3. In the same way, if $X$ is not realcompact, it is
possible to give linear biseparating maps $T: C(X) \ra C(Y)$
which are not continuous (see for instance \cite{ABN1}).

\section{An application to linear isometries between spaces of vector-valued
bounded continuous  functions}

In this section we assume that $E$ and $F$ are ${\Bbb K}$-Banach spaces, and that we are in one of the following situations.

\begin{itemize}
\item {\bf Situation 1.} $\bx = C_b(X, E)$ and $\ay = C_b(Y, F)$. In this case we also assume that both $X$ and $Y$ are realcompact spaces, and that both $E$ and $F$ are infinite-dimensional.
\item {\bf Situation 2.} $\bx= C_b^u(X, E)$ and $\ay = C_b^u(Y, F)$. In this case $X$ and $Y$ will be complete metric spaces.
\end{itemize}

As in the previous section, we also assume that if we are in Situation 1, then $A = C_b (X)$. Otherwise  $A= C_b^u(X)$.

For a ${\Bbb K}$-Banach space $B$, we denote by ${\rm Ext}_B$ the set of extreme points of the closed unit ball of its dual $B'$.

\begin{definition}
{\em Given a Banach space $B$, a continuous linear operator $T: B \ra B$ is said to be a multiplier if every $p \in {\rm Ext}_B$ is an eigenvector for $T'$, i.e. if there is a function $a_T : {\rm Ext}_B \ra {\Bbb K}$ such that $p \circ T = a_T (p) p$ for every $p \in {\rm Ext}_B$.}
\end{definition}

The following characterization of multipliers can be found in \cite[Theorem 3.3]{B}.

\begin{theorem}\label{boundi}
Given a Banach space $B$, a continuous linear  operator $T: B \ra B$ is a multiplier if and only if there is a $\lambda >0$ such that, for every $x \in B$, $Tx$ is contained in every ball which contains $\{ \mu x : \mu \in {\Bbb K}, \va \mu \vb \le \lambda \}$.
\end{theorem}

\begin{definition}
{\em Let $B$ be a Banach space. Given two multipliers $T, S: B \ra B$, we say that $S$ is an adjoint for $T$ if $a_S = \overline{a_T}$. 

The centralizer of $B$ is the set of those multipliers $T: B \ra B$ for which an adjoint exists.}
\end{definition}

When it exists, the adjoint operator for $T$, which must be  unique, will be denoted by $T^*$. On the other hand, the centralizer of $B$ will be denoted by $Z(B)$ (notice that when ${\Bbb K} ={\Bbb R}$, the centralizer of $B$ consists of the set of all multipliers for $B$).

\smallskip

Given $h  \in A$, we define the operator  $M_h : \bx \ra \bx$ as $M_h (f) := hf $ for each  $ f \in \bx$.

With a proof similar to that of \cite[Proposition 4.7 (i)]{B}, we have the following result.

\begin{lemma}\label{subh}
For each $h \in  A$, the operator $M_h$ belongs to  $Z(\bx)$.
\end{lemma}

On the other hand, as it is seen in \cite[Proposition 4.7 (i)]{B}, when ${\Bbb K} ={\Bbb C}$, the adjoint for $M_h$ is $M_{\overline{h}}$. Next we are going to see that the converse of this lemma is also true when $Z(E)$ is trivial.

\begin{lemma}\label{bg}
Suppose that $Z(E)$ is one-dimensional. Given an operator $T \in Z(\bx)$, there exists  $h \in A$ such that $T= M_h$.
\end{lemma}

{\em Proof.} 
Suppose  that $T \in Z(\bx)$. Then there exists a map \[a_T : {\rm Ext}_{\bx} \ra {\Bbb K}\] such that $q \circ T = a_T (q) q$ for every $q \in {\rm Ext}_{\bx}$.

 Fix $x \in X$ and define $T_x : E \ra E$ as $T_x  {\bf e} := (T \widehat{\bf e}) (x)$ for each ${\bf e} \in E$, that is, $T_x = e_x \circ T \circ {\bf i}$, where $e_x : \bx \ra E$ is the evaluation map at $x$, and ${\bf i}: E \ra \bx$ is the natural embedding. We are going to prove that $T_x$ is a multiplier for $E$.

	First, we are going to see that if $p \in {\rm Ext}_E$, then $p \circ e_x \in {\rm Ext}_{\bx}$. Imagine that $p \circ e_x  = \alpha p_1 + (1- \alpha) p_2$, where $p_1 , p_2 $ are points in the closed unit ball of the dual space $\bx'$, and $0< \alpha <1$. We have to prove that $p_1 = p_2 = p \circ e_x $.

Notice that $\bx$ can be expressed as the direct sum of the closed subspaces $E_1 := \{\widehat{\bf e} : {\bf e} \in E \}$ and $E_2 := \{ f \in \bx : f(x) =0 \}$. It is easy to see that if we define, for $i=1,2$,  $q_i : E \ra {\Bbb K}$ in such a way that $q_i ({\bf e}) = p_i (\widehat{\bf e})$, then $q_1$ and $q_2$ belong to the closed unit ball of $E'$. Clearly we have that $p  = \alpha q_1 + (1- \alpha) q_2$, and since $p \in {\rm Ext}_E$, we deduce that $q_1 = q_2 =p$. This implies that $p_1 = p_2 =  p \circ e_x  $ in the subspace $E_1$. But this allows us to claim that, if we assume that $p_1 \neq p_2$, then there exists $f_0 \in E_2$ with $f_0 \notin {\rm Ker} \hspace{.02in} p_1$, and $f_0 \in {\rm Ker} \hspace{.02in} p_2$. As a consequence, since $f_0 (x) =0$,
\[ 0 = (p \circ e_x ) (f_0) = \alpha p_1 (f_0), \] which gives $\alpha =0$, against our hypothesis. As a consequence $p_1 = p_2 = p\circ e_x$, and $p \circ e_x \in {\rm Ext}_{\bx}$.

This implies in particular that $(p \circ e_x ) \circ T = a_T (p \circ e_x) p\circ e_x$, which gives us that, for every $f \in \bx$, \[ p ((Tf)(x) ) = a_T (p \circ e_x) p(f(x)).\]

	Consequently we have that, whenever $p \in {\rm Ext}_E$ and ${\bf e} \in E$, \[p (T_x ({\bf e})) = a_T (p \circ e_x) p({\bf e}).\] In this way, if we define $a_{T_x} :{\rm Ext}_E \ra {\Bbb K}$ as $a_{T_x} (p) := a_{T} (p \circ e_x)$, then we will have that, for every $p \in {\rm Ext}_E$, \[p \circ T_x = a_{T_x} (p) p.\] As a consequence, $T_x$ is a multiplier.

	But notice that working as above we can prove that the operator $e_x \circ T^* \circ {\bf i} : E \ra E$ is also a multiplier. On the other hand it is straightforward to see that it is the adjoint for $T_x$. Consequently $T_x $ belongs to $Z(E)$.

	Now, as $Z(E) = {\Bbb K} \hspace{.02in} {\bf Id}_E$, we have that there exists $a_x \in {\Bbb K}$ such that $T_x = a_x {\bf Id}_E$, and this implies clearly that, for every $p \in {\rm Ext}_E$, $a_{T_x} (p) = a_x$, that is, $a_{T_x}$ is a constant function.

	Thus, given $f \in \bx$, we saw above that, for every $p \in {\rm Ext}_E$, $ p ((Tf)(x) ) = a_T (p \circ e_x) p(f(x))$, that is, \[ p ((Tf)(x) ) = a_x p(f(x)) =p (a_x f(x)).\] This clearly implies that \[(Tf) (x) = a_x f(x),\] because ${\rm Ext}_E$ separates the points of $E$. Since this is true for every $x \in X$, we conclude that, if we define $h: X \ra {\Bbb K}$ as $h(x) := a_x$ for each $x \in X$, then $Tf = hf$ for every $f \in \bx$. Finally, it is clear that, given ${\bf e} \in E$, ${\bf e} \neq 0$, $T \widehat{\bf e} = h \widehat{\bf e}$, which belongs to $\bx$. Then it is easy to prove that $h \in A$, and now we can say that $T = M_h$.
\hfill $\Box$

\bigskip

In the next three results, we assume $B=C_b(Y)$ if we are in Situation 1, and $B=C_b^u(Y)$ if in Situation 2.

The proof of the following  lemma is an adaptation of the one given in \cite[Lemma 4.13 (i)]{B}.

\begin{lemma}\label{dan}
If $T: \bx \ra \ay$ is a surjective linear isometry, then for each  $h \in B$, the map $\widehat{T} M_h $, defined as $(\widehat{T} M_h ) (f) := T^{-1} (h Tf)$ for each  $f \in \bx$, belongs to  $Z(\bx)$.
\end{lemma}

{\em Proof.} 
First notice that $(T^{-1}) ' : \bx ' \ra \ay'$ is a linear surjective isometry, and consequently it maps ${\rm Ext}_{\bx} $ onto ${\rm Ext}_{\ay}$. So, for $p \in {\rm Ext}_{\bx}$, $(T^{-1})' (p) $ belongs to ${\rm Ext}_{\ay}$, and then, by Lemma~\ref{subh}, it is clear that $((T^{-1})' (p)) \circ M_h = a_{M_h} ((T^{-1})' (p)) \cdot (T^{-1})' (p)$. But this means that, for every  $p \in {\rm Ext}_{\bx}$, 
$p \circ (T^{-1}  \circ M_h \circ T) = a_{M_h} (p \circ T^{-1}) \cdot p$, which implies that $T^{-1}  \circ M_h \circ T$ is a multiplier such that $a_{T^{-1}  \circ M_h \circ T} (p) =a_{M_h} (p \circ T^{-1}) $.

So the lemma is proved if ${\Bbb K} = {\Bbb R}$. Now, if ${\Bbb K} = {\Bbb C}$, 
we just have to find an adjoint for $T^{-1}  \circ M_h \circ T$. But notice that if $h \in B$, then $\overline{h}$ also belongs to $B$. We deduce that, in the same way as above, $T^{-1}  \circ M_{\overline{h}} \circ T$ is also a multiplier and, for every $p \in {\rm Ext}_{\bx}$, $
a_{T^{-1}  \circ M_{\overline{h}} \circ T} (p) = a_{M_{\overline{h}}} (p \circ T^{-1}) $. Finally, since the adjoint for $M_h$ is $M_{\overline{h}}$, as we remarked after Lemma~\ref{subh}, we conclude that $T^{-1}  \circ M_{\overline{h}} \circ T$ is the adjoint for $T^{-1}  \circ M_h \circ T$, and we are done.
\hfill $\Box$ 
\bigskip

In the proof of the next theorem,  when our spaces satisfy Situation 1, then $\gamma X$ and $\gamma Y$
are just the Stone-\v{C}ech compactifications of $X$ and $Y$, respectively. 

As for the case when we are in Situation 2, notice that every $f \in C_b^u(X)$ can be extended to a continuous map from $\beta X$ into ${\bf K}$. Also, we can introduce an equivalence relation $\sim$ in $\beta X$ as follows: $x \sim y$ whenever the extensions of every $f \in C_b(X)$ take the same values at $x$ and at $y$. In this way, we obtain the quotient space $\gamma X:= \beta X/ \sim$, which is a compactification of $X$.  Notice that every $f \in C_b^u(X)$ can be extended to  a continuous map from $ \gamma X$ into ${\Bbb K}$.  Then $C_b^u(X)$ can be isometrically embedded in $C(\gamma X)$. Also, we can easily see that the Stone-Weierstrass theorem implies that $C_b^u(X)$ and $C(\gamma X)$ coincide.

On the other hand, a similar process leads to the definition of $\gamma Y$, and this will allow us to identify  $B $ with $C(\gamma Y)$.

\begin{proposition}\label{bono}
Suppose that $Z(E)$ and $Z(F)$ are one-dimensional. If  $T: \bx \ra \ay$ is a surjective linear isometry, then $T$ is  biseparating.
\end{proposition}

{\em Proof.} 
Clearly, it is enough to prove that $T^{-1}$ is  separating, because a similar argument would allow us to conclude that $T$ is also   separating.

First we are going to prove that if there exist two functions $g_1$ and $g_2$ in $\ay$ with $c(g_1) \cap c(g_2) = \emptyset$ and such that $c(T^{-1} g_1) \cap c(T^{-1} g_2) \neq \emptyset$, then we can also find another functions $f_1$ and $f_2$ in $\ay$ with $\cgy c(f_1) \cap \cgy c(f_2) = \emptyset $ and 
$c(T^{-1} f_1) \cap c(T^{-1} f_2) \neq \emptyset$. To this end, assume that $\vc      (T^{-1} g_1) (x_0)  \vd =r_1 >0$ and $\vc (T^{-1} g_2) (x_0) \vd =r_2 >0$ for some $x_0 \in X$. Now let $r_0 := \min \{r_1, r_2\}$, and take $k_1  \in B$ such that $0 \le k_1  \le 1$, $k_1 \equiv 1$ on $ \{ y \in Y  :  \vc g_1 (y) \vd \ge 
r_0 /2 \}$ and $k_1 \equiv 0$ on  $ \{ y \in Y  :  \vc g_1 (y) \vd \le 
r_0 /4 \}$. In the same way, 
take $k_2  \in B$ such that $0 \le k_2 \le 1$, $k_2 \equiv 1$ on $ \{ y \in Y  :  \vc g_2 (y) \vd \ge 
r_0 /2 \}$ and $k_2 \equiv 0$ on  $ \{ y \in Y  :  \vc g_2 (y) \vd \le 
r_0 /4 \}$.  Now define $f_i := k_i g_i$, $i =1,2$. It is clear that $\vc f_i - g_i \vd \le r_0/2$, $i=1,2$. Consequently, as $T^{-1}$ is an isometry, we see that $(T^{-1} f_1) (x_0) \neq 0$ and $(T^{-1} f_2) (x_0) \neq 0$. Also it is easy to prove that $\cgy c(f_1) \cap \cgy c(f_2) = \emptyset$.

Then consider $f_1 , f_2 \in \yf$ such that $\cgy c(f_1) \cap \cgy c(f_2) = \emptyset$. Next take $h \in B$ such that $h \equiv 0$ on $\cgy c(f_2)$ and $h \equiv 1$ on $\cgy c(f_1)$. On the one hand, by Lemmas~\ref{bg} and \ref{dan} we know that $\widehat{T} M_h = M_g$ for some $g \in C( \gamma X)$, that is, \[ (\widehat{T} M_h ) (T^{-1} f_1) = g T^{-1} f_1\] and \[ (\widehat{T} M_h ) (T^{-1} f_2) = g T^{-1} f_2.\]

On the other hand, by definition, 
\[ (\widehat{T} M_h ) (T^{-1} f_1) =  T^{-1} (h f_1) = T^{-1} (f_1)\]
and 
\[ (\widehat{T} M_h ) (T^{-1} f_2) =  T^{-1} (h f_2) = 0.\]

This implies that  $g T^{-1} f_1 =  T^{-1} f_1  $ and $g T^{-1} f_2 = 0$, which gives us that $g=1$ on $c(T^{-1} f_1)$ and $g=0$ on $c(T^{-1} f_2)$. Consequently, by the above reasoning, $T^{-1}$ is  separating, as we wanted to see.
\hfill $\Box$

\begin{definition}
{\em A surjective linear isometry $\txy$ is said to be a 
{\em strong Banach-Stone map} if
there exist a continuous map $J: Y \ra I(E,F)$ and a 
homeomorphism $h:Y \ra X$ such that for every $f \in \bx$ and $y
\in Y$, $(Tf)(y) = (Jy) (f(h(y))$.

A strong Banach-Stone map $T: C_b^u(X, E) \ra C_b^u(Y, F)$ is said to be  {\em uniform} if both $h$ and $h^{-1}$ are uniformly continuous.}
\end{definition}

        Now we are in a position to prove the following
theorem.

\begin{theorem}\label{gol}
Let $E$ and $F$ be Banach spaces such that $Z(E)$ and $Z(F)$ are
one-dimensional. If $T :\xe \ra \yf$ (respectively, $T: C_b^u(X, E) \ra C_b^u(Y, F)$) is a surjective linear isometry,
then it is a strong Banach-Stone map (respectively, a uniform strong Banach-Stone map).
\end{theorem}

{\em Proof.} 
We start with the spaces $\xe$ and $\yf$. Since $T$ is  biseparating, we have that there exist a
homeomorphism $h: Y \ra X$ and a continuous map $J:Y \ra B(E,F)$ such that
for every $f \in \xe$ and $y \in Y$, \[(Tf) (y) = (J y)
(f(h(y)).\] Clearly, we have just to prove that, for each $y \in
Y$, $Jy \in I(E,F)$. Take any $y \in Y$ and ${\bf e} \in E$. Then
$\vc (Jy) (\widehat{\bf e})\vd  = \vc (T
\widehat{\bf e}) (y) \vd$. We are going to see that $ \vc (T
\widehat{\bf e}) (y) \vd = \vc {\bf e} \vd$. Of course, if this is not the case for some $y_0 \in Y$, then $ \vc (T
\widehat{\bf e}) (y_0) \vd < \vc {\bf e} \vd$. Let $C$ an open neighborhood of $y_0 $ in $Y$ such that
\[\sup_{y \in C} \vc (T
\widehat{\bf e}) (y) \vd < \vc {\bf e} \vd.\] Now take $g \in B$ such that $0 \le g \le 1$, $g(Y-C) \equiv 0$, and $g(y_0) =1$. It is clear that there exists $\alpha>0$ such that $ \vc (1+ \alpha g(y)) (T
\widehat{\bf e}) (y)  \vd \le \vc {\bf e} \vd$ for every $y \in Y$. But, on the other hand, it is easy to see that, if $K$ is defined as in the proof of Theorem~\ref{linear}, then
\begin{eqnarray*}
\vc (T^{-1} (1+ \alpha g) (T
\widehat{\bf e}))  (h^{-1} (y_0))  \vd &=& \vc (K h^{-1} (y_0)) ((1+ \alpha g (y_0)) (T \widehat{\bf e}) (y_0)) \vd  \\ & =& (1+ \alpha g (y_0)) \vc (K h^{-1} (y_0)) ( (T \widehat{\bf e}) (y_0)) \vd  \\ &=&  (1+ \alpha g (y_0)) \vc {\bf e} \vd \\ &>& \vc {\bf e} \vd,
\end{eqnarray*}
contradicting the fact that $T$ is an isometry. As a consequence $ \vc (T
\widehat{\bf e}) (y) \vd = \vc {\bf e} \vd$ for every $y \in Y$, and we are done.

Of course, a similar reasoning allows us to prove that every linear surjective isometry $T: C_b^u(X, E) \ra C_b^u(Y, F)$ is a strong Banach-Stone map.
But, in this case, by Corollary~\ref{unif}, we also deduce that it is uniform. 
\hfill $\Box$

\vspace{.15in}

        {\footnotesize {\sc Departamento de Matem\'aticas,
Estadistica y Computaci\'on, Facultad de Ciencias, Universidad
de Cantabria, Avenida de los Castros s.n., E-39071, Santander, Spain.}}

        {\footnotesize {\em E-mail address}: araujoj@@unican.es}

\end{document}